 \documentclass{amsart}

\makeindex

\usepackage[all]{xy}
\usepackage{amssymb}

\usepackage{amsmath}

\usepackage{amsthm}

\usepackage{amscd}

\usepackage{amsfonts}

\usepackage{amssymb}

\newtheorem{theorem}{Theorem}[section]
\newtheorem{lemma}[theorem]{Lemma}
\newtheorem{proposition}[theorem]{Proposition}

\theoremstyle{corollary}
\newtheorem{corollary}[theorem]{Corollary}
\theoremstyle{definition}     
\newtheorem{definition}[theorem]{Definition}

\newtheorem{example}[theorem]{Example}

\theoremstyle{remark}
\newtheorem{remark}[theorem]{Remark}

\numberwithin{equation}{section}

\begin{document}

\title[number of singular points on rational homology projective planes]{The maximum number of singular points on rational homology projective planes}

\author{Dongseon Hwang}
\address{Department of Mathematics, Korea Advanced Institute of Science and Technology, Daejon, Korea}
\email{themiso@kaist.ac.kr}

\author{Jonghae Keum}
\address{School of Mathematics, Korea Institute For Advanced Study, Seoul 130-722, Korea}
\email{jhkeum@kias.re.kr}
\thanks{Research supported by the SRC Program of Korea Science and Engineering Foundation (KOSEF) funded by the Korea government(R11-2007-035-02001-0)}

\subjclass[2000]{Primary 14J17, 14J28}

\date{March 28, 2008.}


\keywords{rational homology projective plane, quotient
singularity, orbifold Bogomolov-Miyaoka-Yau inequality, integral
quadratic form, Enriques surface, symplectic orbifold}

\begin{abstract}
A normal projective complex surface is called a rational homology
projective plane if it has the same Betti numbers with the complex
projective plane $\mathbb{C}\mathbb{P}^2$. It is known that a
rational homology projective plane with quotient singularities has
at most $5$ singular points. So far all known examples have at
most $4$ singular points. In this paper, we prove that a rational
homology projective plane $S$ with quotient singularities such
that $K_S$ is nef has at most $4$ singular points except one case.
The exceptional case comes from Enriques surfaces with a
configuration of 9 smooth rational curves whose Dynkin diagram is
of type $ 3A_1 \oplus 2A_3$.

We also obtain a similar result in the differentiable case and in
the symplectic case under certain assumptions which all hold in
the algebraic case.
\end{abstract}

\maketitle



\section{Introduction}
A normal projective complex surface is called a rational homology
projective plane if it has the same Betti numbers with the complex
projective plane $\mathbb{C}\mathbb{P}^2$. A normal projective
complex surface with quotient singularities is a rational homology
projective plane, if its second Betti number is equal to 1
(\cite{Kollar06}, p. 2). If a rational homology projective plane
is smooth, then it is either $\mathbb{C}\mathbb{P}^2$ or a fake
projective plane, i.e. a smooth projective surface of general type
with $p_g = q = 0$, $K^2=9$.

 Now let $S$ be a rational homology projective plane with quotient singularities. Assume that $S$ is singular.
 L. Brenton constructed such surfaces \cite{Brenton}, and all
examples produced by his method have at most $4$ singular points
\cite{BDP}. On the other hand, from the orbifold
Bogomolov-Miyaoka-Yau inequality (\cite{Sakai}, \cite{Miyaoka},
\cite{Megyesi}), one can derive that $S$ has at most $5$ singular
points, see Corollary \ref{bound5}. However, there has been no
known examples with $5$ singular points. Our main result is :

\begin{theorem}\label{main} Let $S$ be a rational homology projective plane with quotient singularities.
Assume that $K_S$ is nef. Then $S$ has at most $4$ singular points
except the following case:

$S$ has 5 singular points of type $3A_1\oplus 2A_3$, and its
minimal resolution $S'$ is an Enriques surface.
\end{theorem}

An example of the exceptional case is given in Example \ref{ex}.

One of the main ingredients in our proof is the orbifold
Bogomolov-Miyaoka-Yau inequality. This is the reason why we need
the nefness of $K_S$. The case where $-K_S$ is ample has been
recently dealt with by G. B. Belousov  \cite{B}. He has proved
that log-Del Pezzo surfaces of Picard number 1 with quotient
singularities have at most $4$ singular points. Thus, Theorem
\ref{main} holds true without the nefness of $K_S$.

\begin{corollary} The following hold true.
\begin{itemize}
\item [(1)] Rational cohomology projective planes with
quotient singularities have at most $4$ singular points except the
case given in Theorem \ref{main}.
\item [(2)] Integral homology projective planes with quotient
singularities have at most $4$ singular points.
\end{itemize}
\end{corollary}

Here, a rational cohomology projective plane is a normal
projective complex surface having the same rational cohomology
ring with $\mathbb{C}\mathbb{P}^2$. A rational homology projective
plane with quotient singularities is a rational cohomology
projective plane. As regards integral cohomology projective planes
with quotient singularities, D. Bindschadler and L. Brenton
\cite{BB} have proved that they have at most one singular point of
type $E_8$.

The problem of determining the maximum number of singular points
on rational homology projective planes with quotient singularities
is related to the algebraic Montgomery-Yang problem (\cite{MY},
 \cite{Kollar06}).

We remark that if a rational homology projective plane $S$ is
allowed to have rational singularities, then there is no bound for
the number of singular points. In fact, there are rational
homology projective planes with an arbitrary number of rational
singularities. Such examples can be constructed by modifying
Example 5 from \cite{Kollar06}: take a  minimal ruled surface
$X\to \mathbb{P}^1$ with negative section $E$, blow up $m$
distinct fibres into $m$ strings of 3 rational curves $(-2)\,
\textrm{---}\,(-1)\, \textrm{---}\,(-2)$, then contract the proper
transform of $E$ with the $m$ adjacent $(-2)$-curves, and also the
$m$ remaining $(-2)$-curves, to get a rational homology projective
plane with $m+1$ rational singularities.

We now present a brief outline of the proof of Theorem \ref{main}.
Assume that our surface $S$ has $5$ singular points. Then from the
weak version of orbifold Bogomolov-Miyaoka-Yau inequality (see
Theorem \ref{bmy2}) we get one of the following cases for the
$5$-tuple consisting of the orders of local fundamental groups of
singular points:
\begin{displaymath}
\begin{array}{llllll}
(2,2,3,3,3),& (2,2,2,4,4),&&\\
(2,2,2,3,3),& (2,2,2,3,4),& (2,2,2,3,5),& (2,2,2,3,6),\\
(2,2,2,2,q)&{\rm for}\,\, q \geq 2.& &\\
\end{array}
\end{displaymath}
Given its minimal resolution $f : S' \rightarrow S$, the
exceptional curves and the canonical class $K_{S'}$ span a
sublattice $R+\langle K_{S'}\rangle$ of the unimodular lattice
$H^2(S', \mathbb{Z})_{free}:=H^2(S', \mathbb{Z})$/torsion, where
$R$ is the sublattice spanned by the exceptional curves. We note
that rank($R+\langle K_{S'}\rangle)=$ rank($R$) if and only if
$K_S$ is numerically trivial (Lemma \ref{detR}). The list above
gives an infinite list of possible cases for $R$. We reduce this
infinite list for $R$ by using the orbifold Bogomolov-Miyaoka-Yau
inequality (Theorem \ref{bmy}) together with detailed information
about quotient singularities (e.g. Lemmas \ref{taupre}, \ref{V},
\ref{Dp}, \ref{d}, Table \ref{EFK}). Here, we also use the fact
that $|\det(R+\langle K_{S'}\rangle)|$ is a square number if $K_S$
is not numerically trivial (Lemma \ref{detR}). The reduced list
(Propositions \ref{cyc}, \ref{ncyc}) is still an infinite list,
but the infinite part comes from singularities of special type
called singularities of type $T_6$. For each of these cases for
$R$, we then use lattice theoretic arguments to show that, except
the two cases $R=3A_1\oplus 2A_3$ or $4A_1\oplus D_5$, either the
lattice $R$ or $R+\langle K_{S'}\rangle$ cannot be embedded into
the unimodular lattice $H^2(S', \mathbb{Z})_{free}$ (see \S6).
Finally, in \S7 we show that the case $R=3A_1\oplus 2A_3$ is
supported by an example, and the case $R=4A_1\oplus D_5$ can be
ruled out by an argument from the classification theory of
algebraic surfaces and the theory of discriminant quadratic forms.

To prove that the lattice $R$ or $R+\langle K_{S'}\rangle$ cannot
be embedded into the unimodular lattice $H^2(S',
\mathbb{Z})_{free}$, we consider the lattice $M=R+\langle
K_{S'}\rangle$ when it is of the same rank as the unimodular
lattice, and $M=R\oplus R^{\perp}$ otherwise, where $R^{\perp}$ is
the orthogonal complement of $R$ in the unimodular lattice. Then
we use the Local-Global Principle together with computation of
$\epsilon$-invariants (in our case $\epsilon_3$-invariants) to
show that $M$ is not isomorphic over $\mathbb{Q}$ to the
unimodular lattice. The most complicated cases are the cases for
$R$ coming from singularities of type $T_6$. We note that in these
cases rank($R+\langle K_{S'}\rangle)=$ rank($R$), hence we have to
consider $M=R\oplus R^{\perp}$.  We handle this infinite case by
using induction on the rank of $R$ (Lemma \ref{T6notinLd}). There
is an alternative equivalent method: one may compute the
discriminant group of $M$ and proceed to show that this group does
not contain an isotropic subgroup of order the square root of its
order. The latter can be done by showing that the $3$-adic part of
the discriminant group of $M$ does not contain an isotropic
subgroup of order the square root of the order of the $3$-adic
part. We do not give a detailed write-up of this computation. It
takes about the same length of computation as that for
$\epsilon_3$-invariants.

Besides using the theory of algebraic surfaces to analyze the two
cases in \S7, we only use topological facts about algebraic
surfaces and quotient singularities. So we can restate Theorem
\ref{main} in the differentiable case as well as in the symplectic
case under certain assumptions which all hold in the algebraic
case, see \S8.

The first six sections of this paper are as follows. In \S2, we
review the classification theory of cyclic quotient surface
singularities, and prove some properties of Hirzebruch-Jung
continued fractions, which play a key role in reducing the list of
possible cases for $R$. In \S3, we review the orbifold
Bogomolov-Miyaoka-Yau inequality and give some information
regarding the sublattice $R+\langle K_{S'}\rangle$. In \S4-\S5, we
obtain a reduced list for $R$. In \S6, we prove that only two
cases for $R$ may occur.

 Throughout this paper, we work over the field
$\mathbb{C}$ of complex numbers.

\bigskip
{\bf Acknowledgements.} We thank J\'anos Koll\'ar for useful
comments and for the suggestion that the statements for the
differentiable case should be added to our original version. We
also thank Jonathan Wahl for helping us to improve the exposition
of the paper.
\section{Cyclic quotient singularity and $T$-singularity}
In this section, we briefly review the classification theory of
cyclic quotient surface singularities. We also prove some
properties of Hirzebruch-Jung continued fractions, which will be
used later.

\par
Let $p \in Sing(S)$ be a cyclic quotient singularity.  The
irreducible components lying over the point $p$ in its minimal
resolution form a string of smooth rational curves
$\overset{-n_1}{\circ}-\overset{-n_2}\circ-\cdots-\overset{-n_l}\circ$,
and their intersection matrix is given by
\begin{displaymath}
M(-n_1, \ldots, -n_l) = \left( \begin{array}{cccccc}
-n_1 & 1 & 0 & \cdots & \cdots &0 \\
1 & -n_2 & 1 & \cdots & \cdots& 0 \\
0 & 1 & -n_3 & \cdots & \cdots& 0 \\
\vdots & \vdots & \vdots & \ddots & \vdots & \vdots\\
0 & 0& 0 & \cdots  & -n_{l-1}&1  \\
0 & 0& 0 & \cdots &1 & -n_l \\
\end{array} \right)
\end{displaymath}
It is known that the order of the local fundamental group $G_p$ is
equal to the absolute value of the determinant of the matrix
$M(-n_1, \ldots, -n_l)$.
\par

A string of smooth rational curves
$\overset{-n_1}{\circ}-\overset{-n_2}\circ-\cdots-\overset{-n_l}\circ$
is also represented by a continued fraction
 \[
[n_1, n_2, ..., n_l]= n_1 - \dfrac{1}{n_2-\dfrac{1}{\ddots -
\dfrac{1}{n_l}}}
 \]
called the Hirzebruch-Jung continued fraction.

\begin{definition}.
\begin{enumerate}\item For rational numbers $n_1, n_2, ..., n_l$, we define
$$q := |\det(M(-n_1,
\ldots, -n_l))|$$
$$ q_{a_1, a_2, \ldots, a_m} := |\det(M')|$$ where
$M'$ is the $(l-m)\times(l-m)$ matrix obtained by deleting\\
$-n_{a_1}, -n_{a_2}, \ldots, -n_{a_m}$ from $M(-n_1, \ldots,
-n_l)$.  For example,
$$q_1 = |\det(M(-n_2, \ldots, -n_l))| \text{ and } q_{1, l} = |\det(M(-n_2, \ldots, -n_{l-1}))|.$$
\item For convenience, we also define $q_{1, \ldots,
l}=|\det(M(\emptyset))|=1$.\\
Note that $$q_1 q_l = q_{1,l} q + 1,
\quad [n_1, n_2, ..., n_l] = \dfrac{q}{q_1}.$$
\end{enumerate}
\end{definition}

The following fact from linear algebra will be used frequently.

\begin{lemma}[\cite{Megyesi}]\label{linalg}
For rational numbers $n_1, n_2, ..., n_l$, the solution of the
matrix equation
$$\left( \begin{array}{cccccc}
-n_1 & 1 & 0 & \cdots &0 \\
1 & -n_2 & 1 & \cdots & 0 \\
0 & 1 & -n_3 & \cdots & 0 \\
\vdots & \vdots & \vdots & \ddots & \vdots\\
0 & \cdots & \cdots & 1 & -n_l \\
\end{array} \right)  \left( \begin{array}{c}
a_1\\
a_2\\
\vdots\\
a_{l-1}\\
a_l\\
\end{array} \right)
=-\left( \begin{array}{c}
n_1 - 2 + u\\
n_2 - 2\\
\vdots\\
n_{l-1}-2\\
n_l-2+v\\
\end{array} \right)
$$
is given by
$$a_i = 1 - \dfrac{(1-u)|\det(M(-n_{i+1}, \ldots, -n_l))|}{|\det(M(-n_{1}, \ldots, -n_l))|} -\dfrac{(1-v)|\det(M(-n_1, \ldots,
-n_{i-1}))|}{|\det(M(-n_{1}, \ldots, -n_l))|}$$ $$ =1 -
\dfrac{(1-u)q_{1,2, \ldots, i}}{q} -\dfrac{(1-v)q_{i,i+1, \ldots,
l}}{q}$$ for $i = 1,2,...,l$ if $q\neq 0$.
\end{lemma}

\par
There is a special type of quotient singularity, called
$T$-singularity. A quotient singularity which admits a
$\mathbb{Q}$-Gorenstein smoothing is called a singularity of class
$T$.

\begin{definition} [\cite{LW}]
Let $\mathcal{H}$ be the set of all Hirzebruch-Jung continued
fractions $ [n_1, n_2, \ldots, n_l]$,
$$\mathcal{H} = \underset{l}{\bigcup} \{ [n_1, n_2, \ldots, n_l]\mid \textrm{all}\,\,
n_j\,\, \textrm{are integers} \geq 2 \}.$$
\begin{enumerate}
\item A function $\tau : \mathcal{H} \rightarrow \mathcal{H}$
defined by $$\tau([n_1, n_2, \ldots, n_l]) = [2, n_1, n_2, \ldots,
n_{l-1}, n_l + 1]$$ is called a $\tau$-operation. \item A reverse
operation is a function $r : \mathcal{H} \rightarrow \mathcal{H}$
defined by
$$r([n_1, n_2,
\ldots, n_l]) = [n_l, \ldots, n_2, n_1].$$
\end{enumerate}
\end{definition}

\begin{theorem} [\cite{LW}, \cite{KSB}, \cite{Manetti}]\label{manetti} For an integer $d> 0$, let $T_d \subset \mathcal{H}$ be the following set of continued
fractions, or singularities
$$T_d = \Big\{[n_1, n_2, \ldots, n_l]=\dfrac{dn^2}{dna-1}\in \mathcal{H}\mid n,a, \textrm{integers}, n>a>0,\, \gcd(n,a)=1\
\Big\}.$$ Then
\begin{enumerate}
\item $[4] \in T_1$, $[3, 3]\in T_2$, $[3, 2, 3]\in T_3$, and $[3, 2, 2, \ldots, 2, 3]$ $(d$ vertices$) \in T_d$.
\item If $x \in T_d$, then $r(x) \in T_d$ and $\tau(x) \in T_d$.
\item Every element of $T_d$ is obtained by starting with one of the singularities described in $(1)$ and iterating $\tau$-operations and reverse operations.
\item If $[n_1, \ldots, n_l] \in T_d$, then $\sum n_j=3l+2-d$.
\item Every singularity of class $T$  is either a rational double
point or a singular point of class $T_d$ for some $d$.
\end{enumerate}
\end{theorem}

Furthermore, Looijenga and Wahl proved that a cyclic quotient
singularity is of class $T$ if and only if $\dfrac{q_1+q_l+2}{q}$
is an integer. More precisely,

\begin{lemma} [\cite{LW}, Proposition 5.9]\label{Td} Let $[n_1, \ldots, n_l]\in \mathcal{H}$.
\begin{enumerate}
\item  $q_1+q_l+2=2q$ if and only if
$[n_1, \ldots, n_l]$ corresponds to a rational double point of
type $A_l$.
\item  $q_1+q_l+2= q$ if and only if $[n_1, \ldots, n_l] \in T_d$ for some $d$.
\end{enumerate}
\end{lemma}


We will also use the following properties of Hirzebruch-Jung
continued fractions.

\begin{lemma}\label{taupre} The value of the formula $ q_1 +q_l - q$ is preserved
under the $\tau$-operation and the reverse operation, where $l$
denotes the length of the corresponding continued fraction.
\end{lemma}

\begin{proof}
Clearly, $ q_1 +q_l - q$ is preserved under the reverse operation.

 Let $v = [n_1,\ldots, n_l]$. Then $\tau(v)= [2, n_1,\ldots, n_{l-1},
 n_l+1]$.
We use the small letter $q$ for $v$ and the capital letter $Q$ for
$\tau(v)$. We see that
$$\dfrac{Q_1}{Q_{1,l+1}}= [n_l+1,n_{l-1}, \ldots, n_1]= [n_l,n_{l-1}, \ldots, n_1]+1 = \dfrac{q}{q_{l}}+1$$
and $Q_{1,l+1} = q_{l}$, thus $Q_1 = q + q_{l}.$   Similarly,
$$\dfrac{Q_{l+1}}{Q_{1,l+1}} = [2, n_1, \ldots, n_{l-1}] = 2-\dfrac{q_{1,l}}{q_{l}},$$
$$\dfrac{Q_{1,2}}{Q_{1,2,l+1}} = [n_l+1,n_{l-1}, \ldots, n_2] =
\dfrac{q_1}{q_{1,l}}+1,$$ hence $$Q_{l+1} = 2q_l-q_{1,l}\,\,
\textrm{ and }\,\, Q_{1, 2} = q_1 + q_{1,l}.$$ Now we have
$$Q_1 + Q_{l+1} - Q =Q_1 + Q_{l+1} -(2Q_1-Q_{1,2})= q_1 + q_l -
q.$$
\end{proof}

\begin{lemma}\label{V} Assume that $l \geq 3$. Let $V_l=\{[n_1, \ldots, n_l]\in \mathcal{H} \mid -1\leq q_1+q_l - q \leq 1 \}$. Then,
the following hold true:
\begin{enumerate}
\item $[2, n_2, \ldots, n_{l-1}, 2] \notin V_l$, \item If $n_1 \geq
3$ and $n_l \geq 3$, then $[n_1, n_2, \ldots, n_{l-1}, n_l] \notin
V_l$, \item No element of $V_l$ satisfies $\sum_{j=1}^l n_j=3l-4$.
\end{enumerate}
\end{lemma}

\begin{proof}
(1) Suppose that $[2, n_2, \ldots, n_{l-1}, 2] \in V_l$. Then $q_1 + q_l - q \leq 1$.\\
Since $q = n_1q_1 - q_{1,2} = 2q_1-q_{1,2}$,\\
$q_l - q_1 = (q_1 + q_l - q) - q_{1,2} \leq 1- q_{1,2} < 0$.\\
On the other hand, since $q = n_l q_l - q_{l-1,l}= 2q_l -
q_{l-1,l}$,\\ $q_1 - q_l= (q_1 + q_l - q) - q_{l-1,l} \leq 1-
q_{l-1,l} < 0$,\\ which is a contradiction. \par
(2) Suppose that $[n_1, n_2, \ldots, n_{l-1}, n_l] \in V_l$. Then $q_1 + q_l - q  \geq -1$.\\
Thus\\
$q_l - q_1 = (q_1 + q_l - q) + (n_1 - 2) q_1 - q_{1,2} \geq -1+
(n_1 - 2) q_1 - q_{1,2}\geq -1+ q_1 - q_{1,2} > 0$.\\
Here, if $-1+ q_1 - q_{1,2} = 0$, then $n_2=n_3=\cdots =n_l=2$, which violates the condition $n_l\geq 3$.\\
On the other hand,\\
$q_1 - q_l = (q_1 + q_l - q) + (n_l - 2) q_l - q_{l-1,l} \geq
-1+q_l - q_{l-1,l} \geq 0$,\\ a contradiction. \par (3) If
$\sum{n_j} = 3l-4$, then $l \geq 4$. Thus no element of $V_3$
satisfies $\sum{n_j} = 3l-4$. We use induction on $l$. Assume that
$l \geq 4$. Assume also that no element of $V_{l-1}$ satisfies
$\sum_{j=1}^{l-1}{n_j} = 3(l-1)-4$.
 If $v=[n_1, n_2, \ldots, n_l] \in V_l$, then by (1) and (2) either
$n_1 =3$ and $n_l = 2$, or $n_1 =2$ and $n_l = 3$. Thus
$v=\tau(v')$ for some $v'\in \mathcal{H}$. Then by Lemma
\ref{taupre}, $v'\in V_{l-1}$. But if $v$ satisfies $\sum{n_j} =
3l-4$, so does $v'$.
\end{proof}

\section{ The orbifold Bogomolov-Miyaoka-Yau inequality}
Let $S$ be a surface with quotient singularities and $f : S'
\rightarrow S$ be a minimal resolution of $S$.

 It is well-known that quotient singularities are log-terminal
 singularities. Thus one can write $$K_{S'} = f^{*}K_S -
 \sum{D_p}$$ where $D_p = \sum(a_jE_j)$ is an effective $\mathbb{Q}$-divisor supported on $f^{-1}(p)=\cup E_j$ and $0 \leq a_j < 1$.
It implies that
\[K^2_S = K^2_{S'} - \sum_{p \in Sing(S)}{D^2_p}.
\]
We also recall the orbifold Euler characteristic
$$ e_{orb}(S) := e(S) - \sum_{p \in Sing(S)} \Big ( 1-\frac{1}{|G_p|} \Big )$$
where $G_p$ is the local fundamental group of $p$.

The following theorem, called the orbifold Bogomolov-Miyaoka-Yau
inequality, is one of the main ingredients in the proof of our
main theorem.

\begin{theorem}[\cite{Sakai}, \cite{Miyaoka}, \cite{KNS}, \cite{Megyesi}]\label{bmy} Let
$S$ be a normal projective surface with quotient singularities
such that $K_S$ is nef. Then
\[
K_{S}^2 \leq 3e_{orb}(S).
\]
\end{theorem}

We also need the following weaker inequality, which also holds
when $K_S$ is nef.

\begin{theorem}[\cite{KM}]\label{bmy2} Let
$S$ be a normal projective surface with quotient singularities
such that $-K_S$ is nef. Then
\[
0\leq e_{orb}(S).
\]
\end{theorem}

We know that the torsion free part of the second cohomology group,
$$H^2(S', \mathbb{Z})_{free}:=H^2(S', \mathbb{Z})/{\rm torsion},$$
has a lattice structure which is unimodular. For a singular point
$p\in S$, let $R_p$ be the sublattice of $H^2(S',
\mathbb{Z})_{free}$ spanned by the numerical classes of the
components of $f^{-1}(p)$. Let
$$R=\oplus_{p \in Sing(S)} R_p$$
be the sublattice of $H^2(S', \mathbb{Z})_{free}$. We also
consider the sublattice $R+\langle K_{S'}\rangle$ of $H^2(S',
\mathbb{Z})_{free}$ spanned by $R$ and the canonical class
$K_{S'}$. Note that $${\rm rank}(R)\le {\rm rank}(R+\langle
K_{S'}\rangle)\le{\rm rank}(R)+1.$$

\begin{lemma}\label{detR} The following hold true.
\begin{enumerate}
\item ${\rm rank}(R+\langle K_{S'}\rangle)={\rm
rank}(R)$ if and only if $K_S$ is numerically trivial.
\item $\det(R+\langle K_{S'}\rangle)=\det(R)\cdot K_S^2$ if $K_S$ is not numerically trivial.
\item If $S$ is a rational homology projective plane with quotient
singularities, and if $K_S$ is not numerically trivial, then
$R+\langle K_{S'}\rangle$ is a sublattice of finite index in the
unimodular lattice $H^2(S', \mathbb{Z})_{free}$, in particular
$|\det(R+\langle K_{S'}\rangle)|$ is a square number.
\end{enumerate}
\end{lemma}

\begin{proof}
(1) follows from the equality $K_{S'} = f^{*}K_S - \sum{D_p}$.\\
(2) follows from the fact that $\sum{D_p}$ is a
$\mathbb{Q}$-linear combination of generators of $R$, and
$f^{*}K_S$ is orthogonal to $R$.\\
(3) follows from (1).
\end{proof}

 The following corollary is well-known.

\begin{corollary}\label{bound5}
A rational homology projective plane $S$ with quotient
singularities has at most $5$ singular points.
\end{corollary}

\begin{proof} Since $b_2(S)=1$, either $K_S$ is nef or $-K_S$ is ample.
Let $f : S' \rightarrow S$ be a minimal resolution of $S$.
Quotient singularities are rational, so $p_g(S') = q(S') = 0$.
Thus, by the Noether formula, $e(S') + K^2_{S'} = 12$.  Theorem
\ref{bmy} or \ref{bmy2} imply that
\[
0 \leq  e_{orb}(S) = 3-\sum_{p \in Sing(S)} \Big (1 -
\dfrac{1}{|G_p|}\Big ).
\]
Thus $S$ has at most $6$ singular points. Assume that $S$ has
exactly $6$ singular points. Then, $|G_p| = 2$ for all $p \in
Sing(S)$ and $b_2(S') = b_2(S)+6=7$. Thus $K^2_{S'} = 3$ by the
Noether formula. The lattice $R+\langle K_{S'}\rangle$ is of
finite index in $H^2(S', \mathbb{Z})_{free}$. Its discriminant
$\det(R+\langle K_{S'}\rangle)=2^63$ is not a square, so it cannot
be embedded into a unimodular lattice of the same signature, a
contradiction.
\end{proof}

\begin{lemma}\label{listorders}.
Let $S$ be a rational homology projective plane with quotient
singularities. Assume that $S$ has $5$ singular points. Then the
$5$-tuple consisting of the orders of local fundamental groups of
singular points is one of the following:
\begin{displaymath}
\begin{array}{llllll}
(2,2,3,3,3),& (2,2,2,4,4),&&\\
(2,2,2,3,3),& (2,2,2,3,4),& (2,2,2,3,5),& (2,2,2,3,6),\\
(2,2,2,2,q)&{\rm for}\,\, q \geq 2.& &\\
\end{array}
\end{displaymath}
\end{lemma}

\begin{proof}
Theorem \ref{bmy} implies that
\[
0 \leq e_{orb}(S)= -2+\sum_{p \in Sing(S)} \dfrac{1}{|G_p|}
\]
from which we obtain the list.
\end{proof}

The list above gives an infinite list of possible cases for $R$.
In the next two sections we will reduce this infinite list for $R$
by using the orbifold Bogomolov-Miyaoka-Yau inequality (Theorem
\ref{bmy}) together with detailed information about quotient
singularities (e.g. Lemmas \ref{taupre}, \ref{V}, \ref{detR}). The
following two lemmas, useful to calculate $K^2_S$, are also part
of such information.

\begin{lemma}[\cite{LW}, Proposition 5.9 (iii)]\label{Dp}
Let $p$ be a cyclic quotient singular point of $S$. Assume that
$f^{-1}(p)$ has $l$ components $E_1, \ldots, E_l$ with ${E_i}^2 =
-n_i$ forming a string of smooth rational curves
$\overset{-n_1}{\circ}-\overset{-n_2}\circ-\cdots-\overset{-n_l}\circ$.
\begin{enumerate}
\item If $l = 1$, then $D^2_p = -\dfrac{(n_1-2)^2}{n_1}$. \item If
$l \geq 2$, then $D^2_p =2l - \sum n_j + a_1 + a_l= 2l - \sum n_j
+ 2 - \dfrac{q_1 + q_l + 2}{q}$.
\end{enumerate}
\end{lemma}

\begin{lemma}\label{d} Let $p$ be a non-cyclic quotient singular point of type $D_{q, q_1}$
with the dual graph given by $<b;2,1;2,1;q, q_1>$ $($see Table 1
for the notion of dual graph$)$. Let $l$ be the length of the
string $<q,
q_1>=\overset{-n_1}{\circ}-\overset{-n_2}\circ-\cdots-\overset{-n_l}\circ$.
Assume that $l \geq 2$. Then we have the following:
\begin{enumerate}
\item $\det(R_p) =
(-1)^{l+3}4\{(b-1)q - q_1\}$.
\item $a_l =1 - \dfrac{(b-1)q_l -q_{1, l}}{(b-1)q-q_1}$
\item $D_p^2 = 2l -\Sigma n_j + a_l - (b-2)$
\end{enumerate}
\end{lemma}

\begin{proof}
(1) is just a linear algebra computation.\\
(2) Since $E_j.K_{S'} = n_j - 2$ by the adjunction formula, we
have the following matrix equation:

\begin{displaymath}
\left( \begin{array}{cccccccc}
-2 & 0 & 1 & 0 & 0 &0& \cdots & 0 \\
0 & -2 & 1 & 0 & 0 &0& \cdots & 0 \\
1 & 1 & -b & 1 & 0 &0& \cdots &0 \\
0 & 0 & 1 & -n_1 & 1 &0& \cdots & 0 \\
0 & 0 & 0 & 1 & -n_2 & 1 &\cdots & 0 \\
0 & 0 & 0 & 0 & 1 & -n_3 &\cdots & 0 \\
\vdots & \vdots & \vdots & \vdots & \vdots & \vdots & \ddots & \vdots\\
0 & 0 & 0 & 0 & 0 & \cdots  & 1 & -n_l \\
\end{array} \right)
\left( \begin{array}{c}
a_{l+1}\\
a_{l+2}\\
a_0\\
a_1\\
a_2\\
a_3\\
\vdots\\
a_l\\
\end{array} \right)
=-\left( \begin{array}{c}
0\\
0\\
b-2\\
n_1-2\\
n_2-2\\
n_3-2\\
\vdots\\
n_l-2\\
\end{array} \right)
\end{displaymath}

We see that $2a_{l+1} - a_0 = 0 = 2a_{l+2} - a_0$, hence $a_{l+1}
= a_{l+2} = \dfrac{1}{2}a_0$. So the third row can be rewritten by
$-(b-1)a_0 + a_1 = -(b - 2)$. Thus the above matrix equation can
be simplified to the following.

\begin{displaymath}
\left( \begin{array}{cccccc}
-(b-1) & 1 & 0 & \cdots &0 \\
1 & -n_1 & 1 & \cdots & 0 \\
0 & 1 & -n_2 & \cdots & 0 \\
\vdots & \vdots & \vdots & \ddots & \vdots\\
0 & \cdots & \cdots & 1 & -n_l \\
\end{array} \right)
\left( \begin{array}{c}
a_0\\
a_1\\
a_2\\
\vdots\\
a_l\\
\end{array} \right)
=-\left( \begin{array}{c}
(b-1)-2+1\\
n_1-2\\
n_2-2\\
\vdots\\
n_l-2\\
\end{array} \right)
\end{displaymath}

Since $l \geq 2$, by Lemma \ref{linalg}, $$a_l = 1 -
\dfrac{|\det(M(-(b-1), -n_1, -n_2, \ldots,
-n_{l-1}))|}{|\det(M(-(b-1), -n_1, -n_2, \ldots, -n_l))|} = 1 -
\dfrac{(b-1)q_l - q_{1,l}}{(b-1)q-q_1}.$$

From the matrix equation, we observe that
\begin{eqnarray*}
D^2_p &=& - \overset{l}{\underset{j = 1}{\sum}}a_j(n_j-2) - a_0(b-2)\\
      &=& -\overset{l}{\underset{j = 1}{\sum}}(n_j-2) - a_0 + a_1 + a_l - a_0(b-2)\\
      & =& - \overset{l}{\underset{j = 1}{\sum}}{n_j} + 2l + a_l - (b - 2).
\end{eqnarray*}

\end{proof}
\section{Case: $S$ with only cyclic quotient
singularities}

Let $S$ be a rational homology projective plane with quotient
singularities. In this section we consider the case when $S$
admits only cyclic quotient singularities.

By $A_n$, $D_n$, $E_n$ we denote the negative definite root
lattices.

\begin{proposition}\label{cyc} Let $S$ be a rational homology projective plane
with only cyclic quotient singularities. Assume that $K_S$ is nef.
Assume that $S$ has $5$ singular points. Then we get one of the
the following cases for $R=\oplus_{p \in Sing(S)} R_p$:\\
$R=3A_1\oplus 2A_3$, $3A_1 \oplus A_2\oplus \langle-5\rangle$,
$3A_1 \oplus A_2\oplus A_4,$ $3A_1 \oplus 2A_2$, $4A_1 \oplus
A_5$, or $4A_1\oplus [n_1, n_2, \ldots, n_l]$ for any $[n_1, n_2,
\ldots, n_l]\in T_6$.
\end{proposition}

\begin{proof} We will use the orbifold Bogomolov-Miyaoka-Yau
inequality (Theorem \ref{bmy}) together with detailed information
about quotient singularities (e.g. Theorem \ref{manetti}, Lemmas
\ref{Td}, \ref{taupre}, \ref{V}, \ref{Dp}, \ref{d}). Here, we also
use the fact that $|\det(R+\langle K_{S'}\rangle)|$ is a square
number if $K_S$ is not numerically trivial (Lemma \ref{detR}). We
consider each of the cases given in Lemma \ref{listorders}.

\medskip
(1) The case $(2,2,3,3,3)$\\
The lattice $R$ is one of the following:
$$2A_1 \oplus 3\langle-3\rangle,\,\, 2A_1 \oplus A_2\oplus 2\langle-3\rangle,\,\, 2A_1 \oplus 2A_2\oplus \langle-3\rangle,\,\,
2A_1 \oplus 3A_2,$$ and $K_{S'}^2=4,3,2,1$, respectively. Using
Lemma \ref{Dp}, we get $K_S^2=5, \frac{11}{3}, \frac{7}{3}, 1$,
respectively. Thus in each case, $K_S^2\neq 0$, hence $K_S$ is not
numerically trivial. Furthermore, $\det(R+\langle
K_{S'}\rangle)=\det(R)\cdot K_S^2=(-2^23^3)5,\,
(2^23^3)(\frac{11}{3}),\, (-2^23^3)(\frac{7}{3}),\, 2^23^3$,
respectively. None of these discriminants is a square number
modulo $\pm$ sign, so the lattice $R+\langle K_{S'}\rangle$ cannot
be embedded into a unimodular lattice of the same signature, a
contradiction.

\medskip
(2) The case $(2,2,2,4,4)$ \\
Here, the lattice $R$ is one of the following:
$$3A_1 \oplus 2\langle-4\rangle,\,\, 3A_1 \oplus A_3\oplus \langle-4\rangle,\,\, 3A_1 \oplus 2A_3,$$
and $K_{S'}^2=4,2,0$, respectively. In the first two cases, by
Lemma \ref{Dp} we see that $K_S^2\neq 0$, and $\det(R+\langle
K_{S'}\rangle)=\det(R)\cdot K_S^2=(-2^34^2)6,\, (-2^34^2)3$,
respectively. None of these is a square number modulo sign. Thus,
the lattice $R+\langle K_{S'}\rangle$ cannot be embedded into a
unimodular lattice of the same signature, a contradiction.

In the last case, $b_2(S')=10$, $K_{S'}=f^*(K_S)$, and hence by
Noether formula $K^2_S = K^2_{S'} = 0$. In particular, $K_{S}$ is
numerically trivial. This gives the first case for $R$.

\medskip
(3) The case $(2,2,2,3,6)$ \\
The lattice $R$ is one of the following:
$$3A_1 \oplus \langle-3\rangle\oplus \langle-6\rangle,\,\, 3A_1 \oplus A_2\oplus \langle-6\rangle,\,\,
3A_1 \oplus \langle-3\rangle\oplus A_5,\,\,3A_1 \oplus A_2\oplus
A_5,$$ and $K_{S'}^2=4,3,0,-1$, respectively. In the first three
cases, we see that $K_S^2\neq 0$, and $\det(R+\langle
K_{S'}\rangle)=\det(R)\cdot K_S^2$ is not a square number modulo
sign. Thus, the lattice $R+\langle K_{S'}\rangle$ cannot be
embedded into a unimodular lattice of the same signature, a
contradiction.

In the last case $K_{S}^2=K_{S'}^2=-1$, a contradiction.

\medskip
(4) The case $(2,2,2,3,5)$ \\
The lattice $R$ is one of the following:\\ $3A_1 \oplus
\langle-3\rangle\oplus \langle-5\rangle,\,\, 3A_1 \oplus A_2\oplus
\langle-5\rangle,\,\,
3A_1 \oplus -3\oplus [3, 2],\,\, 3A_1 \oplus A_2\oplus [3, 2],$\\
$3A_1 \oplus \langle-3\rangle\oplus A_4,\,\,3A_1 \oplus A_2\oplus
A_4,$\\ and $K_{S'}^2=4,3,3,2,1,0$, respectively. Except the
second and the last case, we see that $K_S^2\neq 0$, and
$\det(R+\langle K_{S'}\rangle)=\det(R)\cdot K_S^2$ is not a square
number modulo sign. Thus, the lattice $R+\langle K_{S'}\rangle$
cannot be embedded into a unimodular lattice of the same
signature.

In the second case, it can be checked that $K_S^2\neq 0$ and
$\det(R+\langle K_{S'}\rangle)=\det(R)\cdot K_S^2$ is a square.
This gives the second case for $R$.

In the last case, $b_2(S')=10$ and $K^2_{S'} = K^2_S = 0$. This
gives the third case.

\medskip
(5) The case $(2,2,2,3,4)$ \\
There are 4 possible cases for $R$. In each case, we see that
$K_S^2\neq 0$, and $\det(R+\langle K_{S'}\rangle)=\det(R)\cdot
K_S^2$ is not a square number modulo sign. Thus, the lattice
$R+\langle K_{S'}\rangle$ cannot be embedded into a unimodular
lattice of the same signature.

\medskip
(6) The case $(2,2,2,3,3)$ \\
There are 3 possible cases for $R$. In each case, we see that
$K_S^2\neq 0$. The absolute value of the discriminant of
$R+\langle K_{S'}\rangle$ is a square only if $R=3A_1 \oplus 2A_2$
with $K_{S'}^2=2$. The latter gives the fourth case.

\medskip
(7) The case $(2,2,2,2,q), \,\,q\geq 2$\\
In this case, we use the nefness of $K_S$. We see that $R = 4A_1
\oplus R_p$ where $|G_p| = q$ and Theorem \ref{bmy} says that $0
\leq K^2_S \leq \dfrac{3}{q}$. Let $l$ be the number of
irreducible components of $f^{-1}(p)$. Note that
$b_2(S')=1+rank(R)=5+l$ and $K^2_{S'} = 5-l$.\\ If $l = 1$, then,
by Lemma \ref{Dp},
$$K^2_S = K^2_{S'} - D^2_p = 4 + \frac{(q-2)^2}{q}
\geq 4, $$
 which is a contradiction.\\ Now assume that $l \geq 2$, and let $[n_1, n_2, \ldots, n_l]$ be
 the Hirzebruch-Jung continued fraction of $R_p$. In
this case, also by Lemma \ref{Dp},

$$K^2_S = K^2_{S'} - D^2_p  = \overset{l}{\underset{j = 1}{\sum}}n_j
-3l + 5 - (a_1 + a_l).$$ So
$$3l -5 + (a_1+a_l) \leq
\overset{l}{\underset{j = 1}{\sum}}n_j \leq \dfrac{3}{q}+ 3l -5 +
(a_1+a_l).$$
 Since $0 \leq a_1+a_l < 2$, we see that $\overset{l}{\underset{j =
1}{\sum}}n_j= 3l-5, 3l-4$, or $3l-3$.

(7-1) Assume that $\overset{l}{\underset{j = 1}{\sum}}n_j= 3l-5$.
Then $a_1=a_l=0$, and $K_S^2=0$. Since $a_1=1-\dfrac{q_1+1}{q}$,
we see that $R_p=A_l$. Then ${\sum}n_j=2l= 3l-5$. Thus, $l=5$ and
$R = 4A_1 \oplus A_5$. This gives the fifth exceptional case.

(7-2) Assume that $\overset{l}{\underset{j = 1}{\sum}}n_j= 3l-3$.
Then by Lemma \ref{linalg}, $$0\leq K^2_S = 2-(a_1+a_l)
=\dfrac{q_1+q_l+2}{q} \leq \dfrac{3}{q}.$$ So $0\leq q_1+q_l+2
\leq 3$, which is impossible.

(7-3) Now assume that $\overset{l}{\underset{j =
1}{\sum}}n_j=3l-4$. First note that ${\sum}n_j=3l-4\geq 2l$, so
$l\geq 4$. By Lemma \ref{linalg},
$$0\leq K^2_S = 1-(a_1+a_l) =\dfrac{q_1+q_l+2}{q}-1 \leq
\dfrac{3}{q}.$$
 Thus
$$q-2\leq q_1+q_l \leq q+1.$$
Hence, by Lemma \ref{V},
 $$q_1 + q_l = q - 2.$$
Then by Lemma \ref{Td}, $[n_1, n_2, \ldots, n_l]\in T_d$ for some
$d$, and by Theorem \ref{manetti}, $d=6$. Furthermore $K_S^2=0$.
This gives the last infinite case for $R$.
\end{proof}

\begin{remark} Except the case $(2,2,2,2,q)$, the
argument above works without the nefness of $K_S$, i.e. works even
in the case when $-K_S$ is ample.
\end{remark}

\begin{remark} Except the two cases $R=3A_1 \oplus A_2\oplus \langle-5\rangle$ and $3A_1 \oplus 2A_2$, we
have shown that ${\rm rank}(R+\langle K_{S'}\rangle)={\rm
rank}(R)$.
\end{remark}
\section{Case: $S$ with a non-cyclic quotient
singularity}

Let $S$ be a rational homology projective plane with quotient
singularities. In this section we consider the case when $S$
admits a non-cyclic quotient singular point.

First we recall Brieskorn's classification of finite subgroups of
$GL(2, \mathbb{C})$ without quasi-reflections \cite{Brieskorn}.
These are generalizations of the famous subgroups of $SL(2,
\mathbb{C})$, i.e. cyclic or binary polyhedral groups. The result
is summarized in Table \ref{eqtable}.

\begin{table}[ht]
\caption{}\label{eqtable}
\renewcommand\arraystretch{1.5}
\noindent\[
\begin{array}{|l|l|l|ll|}
\hline
$Type$&G&|G|&$Dual Graph $ \Gamma_G&\\
\hline
A_{q, q_1}&C_{q, q_1}&q&<q,q_1>&q_1<q,\,\, \gcd(q,q_1)=1\\
\hline
D_{q, q_1}&(Z_{2m}, Z_{2m};D_q, D_q) & 4mq &<b;2,1;2,1;q,q_1> &m=(b-1)q-q_1\,\,$odd$\\
 \hline
D_{q, q_1}&(Z_{4m}, Z_{2m};D_q, C_{2q})&4mq&<b;2,1;2,1;q,q_1> &
m=(b-1)q-q_1\,\,$even$\\
\hline
T_m& (Z_{2m}, Z_{2m};T,T)&24m&<b;2,1;3,2;3,2> & m=6(b-2)+1\\
& &&<b;2,1;3,1;3,1> & m=6(b-2)+5\\
\hline
T_m& (Z_{2m}, Z_{2m};T,D_2)&24m&<b;2,1;3,1;3,2> & m=6(b-2)+3\\
\hline

& &&<b;2,1;3,2;4,3> & m=12(b-2)+1\\
O_m& (Z_{2m}, Z_{2m};O,O)&48m&<b;2,1;3,1;4,3> & m=12(b-2)+5\\
& &&<b;2,1;3,2;4,1> & m=12(b-2)+7\\
& &&<b;2,1;3,1;4,1> & m=12(b-2)+11\\
\hline

& &&<b;2,1;3,2;5,4> & m=30(b-2)+1\\
& &&<b;2,1;3,2;5,3> & m=30(b-2)+7\\
& &&<b;2,1;3,1;5,4> & m=30(b-2)+11\\
I_m& (Z_{2m}, Z_{2m};I,I)&120m&<b;2,1;3,2;5,2> & m=30(b-2)+13\\
& &&<b;2,1;3,1;5,3> & m=30(b-2)+17\\
& &&<b;2,1;3,2;5,1> & m=30(b-2)+19\\
& &&<b;2,1;3,1;5,2> & m=30(b-2)+23\\
& &&<b;2,1;3,1;5,1> & m=30(b-2)+29\\
\hline
\end{array}
\]
\end{table}
 Here we only explain the notation for dual graph.
$$
\begin{array}{lcl}
 <q,q_1> &:=& \text{ the dual graph of the singularity of type }
 \dfrac{1}{q}(1,q_1)
,\\
<b;s_1, t_1; s_2, t_2;s_3, t_3> &:=& \text{ the tree of the
form}\\ & &
\begin{picture}(40,30)
\put(40,25){$<s_2, t_2>$}
 \put(61,15){\line(0,1){6}}
\put(0,5){$<s_1,t_1> -\underset{-b}\circ - <s_3, t_3>$}
\end{picture}
\end{array}
$$

For more information about the table, we refer to the original
paper of Brieskorn\cite{Brieskorn}, Matsuki's
exposition\cite{Matsuki}, or Riemenschneider's work\cite{Rie77}.

\begin{proposition}\label{ncyc}
Let $S$ be a rational homology projective plane with quotient
singularities. Assume that $K_S$ is nef. Assume that $S$ has $5$
singular points including at least one non-cyclic quotient
singular point. Then $R=4A_1\oplus D_5$.
\end{proposition}

\begin{proof} Since  $S$
has a non-cyclic quotient singular point, the possible 5-tuples
are $(2,2,2,2,h).$ In particular, $S$ has only one non-cyclic
quotient singular point, and Theorem \ref{bmy} gives the
inequality

\begin{equation}\label{ncbmy}
0 \leq K^2_S \leq \dfrac{3}{h} \leq \dfrac{3}{8}.
\end{equation}

Let $p\in S$ be the non-cyclic quotient singular point.

\medskip
(1) The case : $p$ is of type $D_{q,q_1}$\\
Let $l$ be the length of the long arm $<q,
q_1>=\overset{-n_1}{\circ}-\overset{-n_2}\circ-\cdots-\overset{-n_l}\circ$
of the dual graph of $f^{-1}(p)$. Then $f^{-1}(p)$ has $l+3$
irreducible components, $b_2(S')=l+8$ and $K_{S'}^2=2-l$.\par If
$l = 1$, then $K^2_{S'} = 1$, thus $1 \leq K^2_S$, a contradiction
to (\ref{ncbmy}). \par
Assume that $l \geq 2$.\\
 By Lemma \ref{d},
 $$K^2_S = K^2_{S'} - D^2_p= \sum{n_j} - 3l + b  - a_l.$$
By (\ref{ncbmy}),  $\sum{n_j} - 3l + b =0$ or 1.

If $\sum{n_j} - 3l + b =0$, then $a_l=0$ and hence by Zariski
lemma (see e.g. \cite{Megyesi}, Lemma 1.3) all components of
$f^{-1}(p)$ are $(-2)$-curves, i.e. $p$ is a rational double
point. Thus $K^2_S =K^2_{S'} =0$. It follows that $l=2$ and $p$ is
of type $D_5$. This gives the case $R=4A_1\oplus D_5$.

If $\sum{n_j} - 3l + b =1$, then
$$K^2_S =1-a_l=
\dfrac{(b-1)q_l - q_{1,l}}{(b-1)q - q_1}\geq \dfrac{1}{(b-1)q -
q_1}=\dfrac{4q}{4mq}\geq\dfrac{8}{4mq}=\dfrac{8}{h},$$ which is a
contradiction to the inequality (\ref{ncbmy}).

\medskip
(2) The case : $p$ is of type $T_m$, $O_m$ or $I_m$\\ By
calculating $K^2_S$ explicitly, we can check that $K^2_S$ does not
satisfy the inequality (\ref{ncbmy}) for every possible case. The
result of exact computation is summarized in Table \ref{EFK}.
\end{proof}

\begin{table}[ht]
\caption{}\label{EFK}
\renewcommand\arraystretch{1.5}
\noindent\[
\begin{array}{|l|l|l|}
\hline
Type & $ Dual Graph $& K^2_S \\
 \hline

T_m &<b;2,1;3,2;3,2>& \frac{6b^2-30b+35}{6b-11}
  \left\{ \begin{array}{ll}
\leq -\frac{1}{7} & $if $ b \leq 3 \\
 \geq \frac{11}{13} & $if $ b \geq 4 \\
\end{array} \right. \\
&<b;2,1;3,1;3,1>&\frac{6b^2-6b-1}{6b-7} \geq \frac{11}{5}\\
\hline
T_m& <b;2,1;3,1;3,2>&\frac{18b^2-54b+41}{18b-27} \geq \frac{5}{9}\\
 \hline

& <b;2,1;3,2;4,3> &\frac{12b^2-72b+94}{12b-23}
  \left\{ \begin{array}{ll}
 \leq -\frac{2}{25} & $if $ b \leq 4 \\
 \geq \frac{34}{37} & $if $ b \geq 5 \\
\end{array} \right. \\
O_m& <b;2,1;3,1;4,3> &\frac{12b^2-48b+46}{12b-19}
  \left\{ \begin{array}{ll}
 =-\frac{2}{5} & $if $ b = 2 \\
 \geq \frac{10}{17} & $if $ b \geq 3 \\
\end{array} \right. \\
& <b;2,1;3,2;4,1> & \frac{12b^2-24b+10}{12b-17} \geq \frac{10}{7}\\
& <b;2,1;3,1;4,1> & \frac{12b^2-14}{12b-13} \geq \frac{34}{11}\\
\hline

& <b;2,1;3,2;5,4> & \frac{30b^2-210b+297}{30b-59}
  \left\{ \begin{array}{ll}
 \leq -\frac{3}{91} & $if $ b \leq 5 \\
 \geq \frac{117}{121} & $if $ b \geq 6 \\
\end{array} \right. \\
& <b;2,1;3,2;5,3> &\frac{30b^2-126b+129}{30b-53}
  \left\{ \begin{array}{ll}
  =-\frac{3}{7} & $if $ b = 2 \\
 \geq \frac{21}{37} & $if $ b \geq 3 \\
\end{array} \right. \\
 & <b;2,1;3,1;5,4> & \frac{30b^2-150b+165}{30b-49}
  \left\{ \begin{array}{ll}
 \leq -\frac{15}{41} & $if $ b \leq 3 \\
 \geq \frac{45}{71} & $if $ b \geq 4 \\
\end{array} \right. \\
I_m& <b;2,1;3,2;5,2> &\frac{30b^2-114b+105}{30b-47}
  \left\{ \begin{array}{ll}
  =-\frac{3}{13} & $if $ b = 2 \\
 \geq \frac{33}{43} & $if $ b \geq 3 \\
\end{array} \right. \\
& <b;2,1;3,1;5,3> &
\frac{30b^2-66b+33}{30b-43} \geq \frac{21}{17}\\
& <b;2,1;3,2;5,1> &\frac{30b^2-30b-15}{30b-41} \geq \frac{45}{19}\\
& <b;2,1;3,1;5,2> &
\frac{30b^2-54b+21}{30b-37} \geq \frac{33}{23}\\
& <b;2,1;3,1;5,1> & \frac{30b^2+30b-63}{30b-31} \geq \frac{117}{29}\\
\hline
\end{array}
\]
\end{table}

\section{Quadratic Forms}
In this section we prove that the cases for $R$ given in
Proposition \ref{cyc} cannot actually occur except the first case
$R=3A_1+2A_3$. We use the Local-Global Principle together with
computation of $\epsilon$-invariants to show that, except the
first case, either the lattice $R$ or $R+\langle K_{S'}\rangle$
cannot be embedded into the unimodular lattice $H^2(S',
\mathbb{Z})_{free}$.
\par

\begin{theorem} \cite{Serre} (Local-Global Principle) Let $f$ and $f'$ be two quadratic forms over $\mathbb{Q}$.
For $f$ and $f'$ to be equivalent over $\mathbb{Q}$ it is
necessary and sufficient that they are equivalent over each
$p$-adic field $\mathbb{Q}_p$ or the field $\mathbb{Q}_{\infty}$
of real numbers.
\end{theorem}

Let $f$ be a quadratic form in $n$ variables over the $p$-adic
field $\mathbb{Q}_p$ such that $f = a_1 {X_1}^2 +  a_2 {X_2}^2 +
\ldots + a_n {X_n}^2$. Define discriminant $d_p(f)$ and
$\epsilon$-invariant $\epsilon_p(f)$ of $f$ as follows:
$$ d_p(f) = a_1 \ldots a_n \in  \mathbb{Q}_p/\mathbb{Q}_p^{*2} $$
$$ \epsilon_p(f) = \underset{i < j}{\prod} (a_i, a_j)_p $$
where $(-, -)_p $ is the Hilbert symbol on $\mathbb{Q}_p$.

Let $f, f'$ be two quadratic forms over the $p$-adic field
$\mathbb{Q}_p$. Then these invariants have the following obvious
properties.
$$ d_p(f\oplus f') = d_p(f)\cdot d_p(f') $$
$$ \epsilon_p(f\oplus f') = \epsilon_p(f) \epsilon_p(f') (d_p(f), d_p(f'))_p$$
We set $\epsilon_p(f) = 1$ if $f$ is a quadratic form in $1$
variable.

\begin{theorem} \cite{Serre} Let $k$ be a $p$-adic field. Then two quadratic forms over
$k$ are equivalent if and only if they have the same rank, the
same discriminant, and the same  $\epsilon$-invariant.
\end{theorem}

Every non-zero element of the $p$-adic field $\mathbb{Q}_p$ can be
written uniquely in the form $p^{\alpha} u$ for some integer
$\alpha$ and some $p$-adic unit  $u$. For any prime number $p$ and
integers $ \alpha$ and $x$ with $1 \leq x < p$, we define
$$  \bar{x} \cdot p^\alpha := \Big \{p^\alpha u \,\,\Big |\,\, u=x + \underset{i \geq 1}{\sum}  a_i p^i \text{ is a
}p \text{-adic unit}\,\, \Big \}.$$

\begin{theorem} \cite{Serre}\label{Serre} (Computation of Hilbert symbol) Let $p > 2$ be a prime number and let $a, b \in \mathbb{Q}_p$.\\
If $a \in \bar{u}\cdot p^\alpha,\,\, b\in \bar{v}\cdot p^\beta$,
then the Hilbert symbol $ (a,b)_p$ can be computed as
$$
(a,b)_p = (-1)^{\alpha \beta \rho (p)}{{u}\choose
\overset{-}{p}}^{\beta} {{v} \choose \overset{-}{p}}^{\alpha}
$$
where ${u \choose \overset{-}{p}}$ denotes the Legendre symbol,
$\rho(p)$ denotes the class modulo $2$ of $\dfrac{p-1}{2}$.
\end{theorem}

\begin{lemma}\label{ortho} Let $L$ be the integral lattice corresponding
 to a Hirzebruch-Jung continued fraction $[n_1, n_2, \ldots, n_l]$ with standard basis $\{e_1, ..., e_l\}$.
Let  $(L\otimes\mathbb{Q}, f)$ be the quadratic form over
$\mathbb{Q}$ defined by $L$. Then we can take an orthogonal basis
$\{v_1, ..., v_l\}$ with ${v_i}^2 = -[n_i, ..., n_1]$ so that the
quadratic form is given by $f = \sum {v_i}^2 {X_i}^2$.
\end{lemma}

\begin{proof}
It is Gram-Schmidt process, essentially.
\end{proof}

\begin{lemma}\label{tauR} Let $L$ be the integral lattice corresponding
 to a Hirzebruch-Jung continued fraction $[n_1, n_2, \ldots, n_l]$. Let $(L\otimes\mathbb{Q}, f_{L})$ be the quadratic from
 with $f_L = \overset{l}{\underset{i = 1}{\sum}} c_i {X_i}^2$
  where $c_i = - [n_i, \ldots, n_1]$ for $i = 1, \ldots, l$.
 Let $(\tau(L)\otimes\mathbb{Q}, f_{\tau(L)})$ be the quadratic form corresponding to $\tau([n_1, n_2, \ldots,
 n_l])$. Then we can choose an orthogonal basis
  $\{v_1, \ldots, v_{l+1}\}$ such that we can write $f_{\tau(L)} = \sum d_i {X_i}^2$
  with $d_i = c_i$ for $i=1, \ldots, l-1$, $d_l = c_l - 1$, and $d_{l+1} =
 -2 - \overset{l}{\underset{j = 1}{\sum}} \dfrac{d_j}{(d_1 d_2 \cdots
 d_j)^2}=-2+\dfrac{q_1+q_{1,l}}{q+q_l}$, where $q=|\det(M(-n_1, -n_2, \ldots, -n_l))|$.

 In particular, if $[n_1, \ldots, n_l]\in T_6$, then the $3$-adic valuation of $ d_{l+1}$ is a positive odd
integer, more precisely, $d_{l+1} \in \bar{2} \cdot 3^{\alpha}$
for a positive odd integer $\alpha$.
\end{lemma}

\begin{proof}
Recall that $\tau([n_1, \ldots, n_l]) = [2, n_1, \ldots, n_{l-1},
n_l+1]$. With respect to a suitable basis $\{ e_1, e_2, \ldots,
e_{l+1}\}$, we can write the corresponding intersection matrix as
follows:
\begin{displaymath}
M_{\tau(L)} = \left( \begin{array}{cccccccc}
-n_1 & 1 & 0 & \cdots & \cdots &0&1 \\
1 & -n_2 & 1 & 0 & \cdots &\cdots &0 \\
0 & 1 & -n_3 & 1 & 0 &\cdots &0 \\
\vdots & \vdots & \ddots & \ddots & \ddots & \ddots & \vdots\\
0 & \cdots & \cdots & 1 & -n_{l-1} & 1 & 0\\
0 & \cdots & \cdots & \cdots & 1 & -n_l-1 & 0\\
1 & 0 & \cdots & \cdots & \cdots & 0 & -2\\
\end{array} \right)
\end{displaymath}
Then, by Gram-Schmidt process, we can write $v_1 = e_1$ and for $i
= 2, \ldots, l+1$,
$$ v_i = e_{i} - \overset{i-1}{\underset{j = 1}{\sum}}
\dfrac{\langle v_j, e_{i} \rangle}{\langle v_j, v_j\rangle} v_j.$$
Then $d_i=v_i^2$ for $i =1, 2, \ldots, l+1$. It is easy to see
that
$$d_i = c_i\,\,\,{\rm for}\,\, i=1, \ldots, l-1,$$
$$d_l = -[n_l+1, n_{l-1}, \ldots, n_1] = -(1 + [n_l, \ldots, n_1])
= -1 + c_l,$$  and
\begin{eqnarray*}
d_{l+1} = e^2_{l+1} - \overset{l}{\underset{j = 1}{\sum}}
\dfrac{{\langle v_j, e_{l+1}\rangle }^2}{\langle v_j, v_j\rangle}
        = -2 - \overset{l}{\underset{j = 1}{\sum}} \dfrac{d_j}{(d_1 \cdots d_j)^2}.
\end{eqnarray*}
Write $c_j = - \dfrac{y_j}{y_{j-1}}$ where $y_0=1$ and
$y_j=|\det(M(-n_1, -n_2, \ldots, -n_j))|$. Clearly $y_l=q$ and
$y_{l-1}=q_l$. Note that $$\dfrac{d_j}{(d_1 \cdots
d_j)^2}=\dfrac{c_j}{(c_1 \cdots c_j)^2} = -\dfrac{1}{y_{j-1}
y_j}$$ for $j = 1, \ldots, l-1$. Claim that
$$\overset{k}{\underset{j = 1}{\sum}}\dfrac{1}{y_{j-1} y_j} =
\dfrac{|\det(M(-n_2, -n_3, \ldots, -n_k))|}{y_k}.$$ We prove the
claim by using induction. If $k = 2$, then
$$\dfrac{1}{y_0y_1} + \dfrac{1}{y_1y_2} = \dfrac{y_2+1}{y_1y_2} =
\dfrac{n_2}{y_2}.$$ Now assume that the claim holds for $k < m$.
Then
\begin{eqnarray*}
\overset{m}{\underset{j = 1}{\sum}}\dfrac{1}{y_{j-1} y_j} & = &
\overset{m-1}{\underset{j = 1}{\sum}}\dfrac{1}{y_{j-1} y_j} +
\dfrac{1}{y_{m-1}y_m}\\
& = &\dfrac{|\det(M(-n_2, -n_3,\ldots, -n_{m-1}))| y_m+1}{y_{m-1}y_m}\\
& = &\dfrac{|\det(M(-n_2, -n_3, \ldots, -n_m))|}{y_m},
\end{eqnarray*}
which proves the claim. Thus
\begin{eqnarray*}
 d_{l+1} & = &-2 - \overset{l-1}{\underset{j =
1}{\sum}}\dfrac{c_j}{(c_1 \cdots c_j)^2} - \dfrac{1}{(c_1 \cdots c_{l-1})^2d_l}\\
& = &-2  + \overset{l-1}{\underset{j =
1}{\sum}}\dfrac{1}{y_{j-1} y_j}- \dfrac{1}{y_{l-1}^2(c_{l}-1)} \\
& = &-2+\dfrac{q_{1,l}}{q_{l}}+\dfrac{1}{q_l(q+q_l)}=-2+
\dfrac{q_1+q_{1,l}}{q+q_l}.
\end{eqnarray*}
Now assume that $[n_1, \ldots, n_l]\in T_6$.  Then
$$q=6n^2,\,\, q_1=6na-1,\,\, q_l=6nb-1$$ for some integers $ n,a,b$
with $n>a>0,\, \gcd(n,a)=1,\, a+b=n$.\\ Since $q_{1,l}q=q_1q_l-1$,
we see that $q_{1,l}=6ab-1$. Using these we get
\begin{eqnarray*}
 d_{l+1}&=& \dfrac{-6(n+b)^2}{6n^2+6nb-1}.
\end{eqnarray*}
Now it is easy to see that $d_{l+1} \in \bar{2} \cdot 3^{\alpha}$
for a positive odd integer $\alpha$.
\end{proof}

\begin{lemma}\label{Lm} \begin{enumerate} \item Let $I_{1,m}: = \langle1\rangle \oplus m\langle-1\rangle$ be the odd unimodular lattice of signature $(1,m)$. Then
$\epsilon_p(I_{1,m}) = 1$ for all $p > 2$.
\item Let $II_{1,8m+1}: = H\oplus mE_8$ be the
even unimodular lattice of signature $(1,8m+1)$, where $H$ is the
even unimodular lattice of signature $(1,1)$, and $E_8$ the even
unimodular lattice of signature $(0,8)$. Then
$\epsilon_3(II_{1,8m+1}) = 1$.
\end{enumerate}
\end{lemma}

\begin{proof} (1) follows from a direct calculation.

(2) It is easy to see that $\epsilon_3(H) = 1$. By a suitable
change of basis, we can write the quadratic form of
$E_8\otimes\mathbb{Q}$ as follows:
$$f = -2X^2_1 - \dfrac{3}{2}X^2_2 - \dfrac{4}{3}X^2_3 -
\dfrac{5}{4} X^2_4 - \dfrac{6}{5}X^2_5 - \dfrac{7}{6}X^2_6-
\dfrac{8}{7}X^2_7 - \dfrac{1}{8}X^2_8.$$
A direct calculation
shows that $\epsilon_3(E_8) = 1$.
 Hence
$$\epsilon_3(H\oplus E_8) =
\epsilon_3(H)\epsilon_3(E_8)(d(H), d(E_8))_3 = 1.$$ Now, use
induction.
\end{proof}

\begin{lemma}\label{T6notinLd} Let $l=m-4$ be an integer $\geq 6$, and  $R_p$ be the lattice of rank $l$ corresponding to
a singularity $p$ of class $T_6$. Then the negative definite
lattice $N: = 4A_1 \oplus R_p$ of rank $m$ cannot be embedded into
the lattice $I_{1,m}$.
\end{lemma}

\begin{proof}
Assume that $N$ is embedded to $I_{1,m}$. Let $N^\perp$ be the
orthogonal complement of $N$ in $I_{1,m}$. Then $(N \oplus
N^\perp) \otimes \mathbb{Q}_3 \cong I_{1,m} \otimes \mathbb{Q}_3.$
Thus by Lemma \ref{Lm},
$$ \epsilon_3(N \oplus N^\perp) = \epsilon_3(I_{1,m})=1.$$
To get a contradiction, we will show that $ \epsilon_3(N \oplus
N^\perp) = -1.$ Note that $\det(N)=(-1)^{l}2^46n^2$ and
$\det(N^\perp)=6n'^2$ for some $n$, $n'$. Hence by Theorem
\ref{Serre}
$$ (d_3(N), d_3(N^\perp))_3 = ((-1)^{l}6, 6)_3=(-1)^{l+1}.$$
It is easy to see that $ \epsilon_3(N)= \epsilon_3(R_p)$. Thus
\begin{eqnarray}
 \epsilon_3(N \oplus N^\perp)
  =  \epsilon_3(N) \epsilon_3(N^\perp) (d_3(N), d_3(N^\perp))_3\nonumber
 = (-1)^{l+1}\epsilon_3(R_p) \nonumber
\end{eqnarray}
It is enough to show that $$\epsilon_3(R_p)=(-1)^{l}.$$ To do this
we use induction on $l$. \par If  $l=6$, then $R_p$ corresponds to
the Dynkin diagram
$\overset{-3}{\circ}-\overset{-2}\circ-\overset{-2}\circ-\overset{-2}\circ-\overset{-2}\circ-\overset{-3}\circ$,
and by Lemma \ref{ortho} the quadratic form
$(R_p\otimes\mathbb{Q}, f)$ over $\mathbb{Q}$ is given by
$$f = -3X^2_1 - \dfrac{5}{3}X^2_2 - \dfrac{7}{5}X^2_3 -
\dfrac{9}{7} X^2_4 - \dfrac{11}{9}X^2_5 - \dfrac{24}{11}X^2_6.$$ A
direct calculation shows that $\epsilon_3(R_p) = 1$.
\par
It is clear that the epsilon invariant does not change under a
reverse operation. Since the $\tau$-operation increases
rank$(R_p)$ by 1, it is sufficient to show that
$$\epsilon_3(R_p) \epsilon_3(\tau(R_p)) = -1.$$
 By Lemma \ref{tauR} and notation there,
$$\epsilon_3(R_p) \epsilon_3(\tau(R_p)) = (c_l d_l d_{l+1}, c_1 \cdots c_{l-1})_3(d_{l+1},
d_l)_3.$$ Recall that $[n_1, \ldots, n_l]=\dfrac{q}{q_1}=
\dfrac{6n^2}{6na-1}$,  $c_l = -[n_l, \ldots, n_1]=-\dfrac{q}{q_l}=
-\dfrac{6n^2}{6nb-1}$,  for some integers $n>a>0$, $n>b>0$ with
$\gcd(n,a) =\gcd(n,b) = 1$. It implies that $q=6n^2$,
$q_1=6na-1\equiv 2 $ mod $ 3$, $q_l=6nb-1\equiv 2 $ mod $ 3$, and
$c_l \in \bar{2} \cdot 3^\alpha$ for some odd integer $\alpha>0$.
Note that
$$c_1 \cdots c_{l-1} = (-1)^{l-1} q_l.$$
\\
Case 1: $l$ is odd.\\
In this case $c_1 \cdots c_{l-1} \in \bar{2} \cdot 3^0$ and $d_l
\in \bar{2} \cdot 3^0$. Thus $(c_l d_l d_{l+1}, c_1 \cdots
c_{l-1})_3 = 1$
and $(d_{l+1}, d_l)_3 = -1$. Hence $\epsilon_3(R_p) \epsilon_3(\tau(R_p))= -1$, as desired.\\
\\
Case 2: $l$ is even.\\
In this case $c_1 \cdots c_{l-1} \in \bar{1} \cdot 3^0$ and $d_l
\in \bar{2} \cdot 3^0$. Thus $(c_l d_l d_{l+1}, c_1 \cdots
c_{l-1})_3 = 1$
and $(d_{l+1}, d_l)_3 = -1$. Hence $\epsilon_3(R_p) \epsilon_3(\tau(R_p)) = -1$, as desired.\\
This completes the proof.
\end{proof}

\begin{lemma}\label{E8} Let
$R=3A_1 \oplus A_2 \oplus A_4$, or $4A_1 \oplus A_5.$\\ Then the
lattice $R$ can be embedded into neither the lattice $I_{1,9}$ nor
$II_{1,9}$.\\
In particular, neither the case $R=3A_1 \oplus A_2 \oplus A_4$ nor
$4A_1 \oplus A_5$  in Proposition \ref{cyc} occurs.
\end{lemma}

\begin{proof}
Suppose that $R$ can be embedded into $L:=I_{1,9}$ or $II_{1,9}$.
 Let $R^\perp$ be the
orthogonal complement of $R$ in $L$. By Lemma \ref{Lm}, it
suffices to show that $\epsilon_3(R\oplus R^\perp)=-1$.\\

Case 1. $R =3A_1\oplus A_2 \oplus A_4$.\\
Since $d(R) = -2^3 \cdot 3 \cdot 5$, we see that $d(R^\perp) =
30$. By a direct calculation, it is easy to see that
$\epsilon_3(R) = -1$, so
$$
\epsilon_3(R\oplus R^\perp)= \epsilon_3(R)
\epsilon_3(R^\perp)(d_3(R), d_3(R^\perp))_3 = -1.$$

Case 2. $R =4A_1\oplus A_5$.\\
Similar to Case 1. Since $d(R) = -2^4\cdot 6$, we see that
$d_3(R^\perp) = 6$. A direct calculation shows that $\epsilon_3(R)
= -1$, so $\epsilon_3(R\oplus R^\perp)=-1$.
\end{proof}

\begin{corollary}\label{enriques} There is no Enriques surface with a configuration of 9 smooth rational curves
whose Dynkin diagram is of type $3A_1 \oplus A_2 \oplus A_4$ or $
4A_1 \oplus A_5.$
\end{corollary}

\begin{proof}
The second cohomology group, modulo torsion, of any Enriques
surface has a lattice structure isomorphic to $II_{1,9}=H \oplus
E_8$.
\end{proof}

\begin{lemma}\label{twocases} The two cases
$R=3A_1 \oplus A_2\oplus \langle-5\rangle$ and $R=3A_1 \oplus
2A_2$ in Proposition \ref{cyc} do not occur.
\end{lemma}

\begin{proof} It suffices to show that the lattice $R=3A_1 \oplus A_2\oplus \langle-5\rangle$ (resp. $3A_1 \oplus
2A_2$ cannot be embedded into the unimodular lattice $H^2(S',
\mathbb{Z})_{free}$ which is isomorphic to the lattice $I_{1,6}$
(resp. $I_{1,7}$). Note that $$(R+\langle
K_{S'}\rangle)\otimes\mathbb{Q}\cong (R+\langle
f^*K_{S}\rangle)\otimes\mathbb{Q}.$$ Thus $$\epsilon_3(R+\langle
K_{S'}\rangle)=\epsilon_3(R+\langle f^*K_{S}\rangle).$$ In case
$R=3A_1 \oplus A_2\oplus \langle-5\rangle$, it can be checked that
$\epsilon_3(R+\langle f^*K_{S}\rangle)=-1$, so by Lemma \ref{Lm}
the lattice
 $R+\langle K_{S'}\rangle$ cannot be embedded into $I_{1,6}=
\langle1\rangle \oplus 6\langle-1\rangle$.

Similarly, in case $R=3A_1 \oplus 2A_2$, it can be checked that
$\epsilon_3(R+\langle K_{S'}\rangle)=-1$, so  by Lemma \ref{Lm}
the lattice $R+\langle K_{S'}\rangle$ cannot be embedded into
$I_{1,7}= \langle1\rangle \oplus 7\langle-1\rangle$.
\end{proof}

Now by Corollary \ref{bound5}, Lemmas \ref{T6notinLd}, \ref{E8},
\ref{twocases}, we can combine Propositions \ref{cyc} and
\ref{ncyc}  into the following form:

\begin{proposition}\label{mainprop} Let $S$ be a rational homology projective plane
with quotient singularities. Assume that $K_S$ is nef. Then $S$
has at most $4$ singular points except the following two cases:

$S$ has $5$ singular points of type $3A_1\oplus 2A_3$  or
$4A_1\oplus D_5$.
\end{proposition}

\begin{proposition}\label{mainprop2} In either case  $R=3A_1\oplus 2A_3$ or $4A_1\oplus D_5$,
$S'$ is an Enriques surface.
\end{proposition}

\begin{proof} In either case, we have shown in the proof of Propositions \ref{cyc} and
\ref{ncyc} that $K_S$ is numerically trivial. Since $S$ has only
rational double points, $K_{S'}=f^*K_S$, hence $K_{S'}$ is
numerically trivial. We know that $p_g(S')=q(S')=0$. Thus by the
classification theory of algebraic surfaces $S'$ is an Enriques
surface.
\end{proof}
\section{Enriques surfaces}

In this section we show that the case $R=3A_1\oplus 2A_3$ is
supported by an example, and the case $R=4A_1\oplus D_5$ can be
ruled out by an argument from the classification theory of
algebraic geometry and the theory of discriminant quadratic forms.

Let $L$ be a non-degenerate even lattice.  The bilinear form of
$L$ determines a canonical embedding $L \subset L^{\ast}= Hom(L,
\mathbb{Z})$. The factor group $L^{\ast}/L$, which is denoted by
$disc(L)$, is an abelian group of order $|\det(L)|$. We denote by
$l(L)$ the number of minimal generators of $disc(L)$. We extend
the bilinear form on $L$ to the one on $L^{\ast}$, taking value in
$\mathbb{Q}$, and define
$$q_{L} : disc(L) \to \mathbb{Q}/2\mathbb{Z},\quad   q_{L}(x+L)\,
=\, \langle x,x \rangle +\,2\mathbb{Z}\,\   (x\in L^{\ast}).$$ We
call $q_{L}$ the {\it discriminant quadratic form} of $L$. A
subgroup $A$ of $disc(L)$ is said to be isotropic if $q_{L}$ takes
value identically 0 on $A$.

For a non-degenerate odd lattice, its discriminant quadratic form
can be defined similarly.

Let $L$ be a sublattice of a lattice $M$. The lattice $L$ is said
to be primitive if $M/L$ is torsion free. The minimal primitive
sublattice of $M$ containing $L$ is called the primitive closure
of $L$, and is denoted by $\bar{L}$. The orthogonal complement of
$L$ in $M$ is denoted by $L^{\perp}_M$, or simply by $L^{\perp}$.
The following is well known (see e.g. \cite{N}).

\begin{lemma}\label{disgr} Let $L$ be a non-degenerate even lattice.
\begin{enumerate}
\item If an even lattice $M$ is an over-lattice of $L$, i.e. $M$ has the
same rank as $L$ and contains $L$, then the group $A:=M/L$ is an
isotropic subgroup of $disc(L)$, and $disc(M)\cong A^{\perp}/A$.
\item Conversely, every
isotropic subgroup $A$ of $disc(L)$ defines a unique over-lattice
$M\subset L^{\ast}$ with $disc(M)\cong A^{\perp}/A$.
\item If $L$ is primitive in a unimodular even lattice, then $$(disc(L^{\perp}), q_{L^{\perp}})\cong
(disc(L), -q_L).$$
\end{enumerate}
\end{lemma}

\begin{proposition}\label{enriq} There is no Enriques surface with a configuration of 9 smooth rational curves
whose Dynkin diagram is of type $ 4A_1 \oplus D_5.$
\end{proposition}

\begin{proof}
Suppose that there is such an Enriques surface $W$. The
N\'eron-Severi group modulo torsion, $H^2(W,
\mathbb{Z})_{free}:=H^2(W, \mathbb{Z})$/torsion, has a lattice
structure isomorphic to $H \oplus E_8$. Here, the torsion is
generated by the canonical class $K_W$. Let $R =4A_1 \oplus D_5$
be the sublattice of $H^2(W, \mathbb{Z})_{free}$ generated by the
9 smooth rational curves on $W$. Let $E_1, E_2, E_3, E_4$ be the
smooth rational curves corresponding to $4A_1$. Note that
$$disc(R) = \Big(\overset{4}{\underset{i = 1}{\oplus}}
(\mathbb{Z}/2) \langle e_i \rangle \Big) \oplus \Big (
(\mathbb{Z}/4)\langle v \rangle \Big ),$$ where $\langle \cdot
\rangle $ is the generator of the group, e.g.
$e_i=\dfrac{[E_i]}{2}$. The quadratic form on $disc(D_5) \cong
(\mathbb{Z}/4)\langle v \rangle $ is given by $v^2=-\dfrac{5}{4}$.
Since rank$R^{\perp}=1$, $disc(R^\perp)$ is a cyclic group. Hence
by Lemma \ref{disgr}, we see that $disc(\bar{R})\cong
-disc(R^\perp)$ is a cyclic group and $\bar{R}/R$ is an isotropic
subgroup of $disc(R)$. Since $l(R)=5$, this is possible only if
$disc(\bar{R}) \cong \mathbb{Z}/4$ and $\bar{R}/R =
\overset{2}{\underset{i = 1}{\oplus}} (\mathbb{Z}/2)$. Finding
generators of $\bar{R}/R$, we see that it is generated by two
elements $e_1 + e_2+e_3+ e_4$ and $e_i + e_j+2f$ for some $i\neq
j$. In any case, $e_1 + e_2+e_3+ e_4 \in \bar{R}/R$. This means
that $E_1 + E_2+E_3+ E_4$ is divisible by 2 in $H^2(W,
\mathbb{Z})_{free}$, i.e. either $E_1 + E_2+E_3+ E_4$ or $E_1 +
E_2+E_3+ E_4+K_W$ is divisible by 2 in $H^2(W, \mathbb{Z})=$
Pic$(W)$.  Let $X$ be the algebraic K3 cover of $W$. Then it
follows that the pre-images in $X$ of the 4 curves $E_1, E_2, E_3,
E_4$ are 8 smooth rational curves whose sum is divisible by 2 in
Pic$(X)$. Let $X\to X'$ be the contraction of these 8 curves. Note
that away from these singular points, $X'$ contains 10 smooth
rational curves whose Dynkin diagram is of type $2D_5$. Then there
is a double cover $Y$ of $X'$ branched exactly along the 8
singular points. The surface $Y$ is an algebraic K3 surface (cf.
\cite{KZ}, Theorem 1 and 2). Then $Y$ contains 20 smooth rational
curves whose Dynkin diagram is of type $4D_5.$ This implies that
$Y$ has Picard number $\ge 21$, which is impossible.
\end{proof}

The following example was mentioned in Theorem \ref{main}.

\begin{example}\label{ex} There is an Enriques surface with a configuration of 9 smooth rational curves
whose Dynkin diagram is of type $ 3A_1 \oplus 2A_3.$ See Example
III, \cite{Kon}. This Enriques surface has an elliptic fibration
with 2 double fibres of type $I_4$, 2 fibres of type $I_2$, and a
special $2$-section intersecting only one component in each fibre.

Let $S$ be a rational homology projective plane with 5
singularities of type $3A_1\oplus 2A_3$. Then $S$ is not an
integral homology projective plane, because $H_1(S,
\mathbb{Z})\cong\mathbb{Z}/2\mathbb{Z}\neq 0$. But $S$ and
$\mathbb{C}\mathbb{P}^2$ have the isomorphic rational cohomology
ring, although $H^2(S, \mathbb{Q})$ does not contain an element of
self-intersection 1.
\end{example}

Now  Theorem \ref{main} follows from Propositions \ref{mainprop},
\ref{mainprop2}, \ref{enriq} and Example \ref{ex}.

\section{The differentiable case}
Let $M$ be a smooth, compact 4-manifold whose boundary components
are spherical, that is, they are links. One can then attach cones
to each boundary component to get a 4-dimensional orbifold $S$. As
in the algebraic case, there is a minimal resolution $f : S'
\rightarrow S$, where $S'$ is a smooth, compact 4-manifold without
boundary.\\
To each singular point $p\in S$ (the vertex of each cone), we
assign a uniquely defined class $D_p = \sum(a_jE_j)\in H^2(S',
\mathbb{Q})$ such that $D_p\cdot E_i = 2+E_i^2$ for each component
$E_i$ of $f^{-1}(p)$.\\
We always assume that $S$ and $S'$ satisfy the following two
conditions:
\begin{itemize} \item [(1)] $S$ is a
$\mathbb{Q}$-homology $\mathbb{C}\mathbb{P}^2$, i.e. $H^1(S,
\mathbb{Q})=0$ and $H^2(S, \mathbb{Q})\cong\mathbb{Q}$.
\item [(2)] The intersection form on $H^2(S', \mathbb{Q})$ is
indefinite, and is negative definite on the subspace generated by
the classes of the exceptional curves of $f$.
\end{itemize} If there is a class $K_{S'}\in H^2(S', \mathbb{Q})$ satisfying  both the
N\"other formula $$ K^2_{S'}=10-b_2(S')$$
 and the adjunction formula $$K_{S'}\cdot E +E^2=-2$$
 for each exceptional curve $E$ of $f : S'
\rightarrow S$, we call it a {\it formal canonical class} of $S'$.

\begin{theorem}\label{diff} Let $M$, $S$, and $S'$ be the same as above satisfying the conditions $(1)$ and $(2)$. Assume that
 $S'$ admits a formal canonical class $K_{S'}$. Assume further that
 \[ K^2_{S'} - \sum_{p \in Sing(S)}{D^2_p}\leq 3e_{orb}(S).
\]
Then $M$ has at most $4$ boundary components except the following
two cases:

$M$ has 5 boundary components of type $3A_1+2A_3$ or $4A_1+D_5$.
\end{theorem}

Note that the assumptions in Theorem \ref{diff} all hold for
algebraic $\mathbb{Q}$-homology projective planes  with quotient
singularities such that the canonical divisor is nef.

\begin{proof} In our proof up to Proposition \ref{mainprop} for the algebraic orbifold case, the canonical class
$K_S$ appears several times, but can be replaced by $f^*K_S$.
Given a formal canonical class $K_{S'}$ in the differentiable
case, the class $K_{S'}+\sum D_p\in H^2(S', \mathbb{Q})$ plays
exactly the same role as $f^*K_S$. The words ``$K_S$ is
numerically trivial" is now replaced by ``$K_{S'}=-\sum D_p \in
H^2(S', \mathbb{Q})$", or by ``$K_{S'}\in R\otimes\mathbb{Q}$".
\end{proof}

\begin{theorem}\label{diff2} Let $M$, $S$, and $S'$ be the same as above satisfying the conditions $(1)$ and $(2)$. Assume that
 $S'$ admits a formal canonical class $K_{S'}$. Assume further that
\[
0 \leq e_{orb}(S).
\]
Then $M$ has at most $5$ boundary components. The bound is sharp.
\end{theorem}

The assumptions in Theorem \ref{diff2} all hold for algebraic
$\mathbb{Q}$-homology projective planes  with quotient
singularities.

If $S$ is a symplectic orbifold, then $S'$ is a symplectic
manifold and the symplectic canonical class $K_{S'}$ gives a
formal canonical class.

\begin{corollary}\label{symp} Let $M$, $S$, and $S'$ be the same as above satisfying the conditions $(1)$ and $(2)$. Assume that
 $S$ is a symplectic orbifold. Assume further that
 \[ K^2_{S'} - \sum_{p \in Sing(S)}{D^2_p}\leq 3e_{orb}(S).
\]
Then $M$ has at most $4$ boundary components except the following
two cases:

$M$ has 5 boundary components of type $3A_1+2A_3$ or $4A_1+D_5$.
\end{corollary}

\begin{corollary}\label{symp2} Let $M$, $S$, and $S'$ be the same as above satisfying the conditions $(1)$ and $(2)$. Assume that
$S$ is a symplectic orbifold. Assume further that
\[
0 \leq e_{orb}(S).
\]
Then $M$ has at most $5$ boundary components. The bound is sharp.
\end{corollary}

\begin{remark}  In the differentiable
case, if a formal canonical class $K_{S'}$ is given, then a {\it
formal canonical class} of $S$ can be defined as the class
$K_{S'}+\sum D_p\in H^2(S', \mathbb{Q})$.
\end{remark}
\bibliographystyle{amsplain}

\begin{thebibliography}{99}


\bibitem {B} G. N. Belousov, \textit{Del-Pezzo surfaces with log-terminal singularities},
Mat. Zametki \textbf{83} (2008), no. 2, 170-180

\bibitem {BB} D. Bindschadler and L. Brenton, \textit{On singular 4-manifolds of the homology type of $CP^2$},
 J. Math. Kyoto Univ. \textbf{24} (1984), 67-81

\bibitem {BDP} L. Brenton, D. Drucker and G. C. E. Prins, \textit{Graph theoretic techniques in algebraic geometry II:
construction of singular complex surfaces of the rational cohomology type of $CP^2$},
Comment. Math. Helvetici \textbf{56} (1981), 39-58

\bibitem {Brenton} L. Brenton, \textit{Some examples of singular compact analytic surfaces
wich are homotopty equivalent to the complex projective plane},
Topology \textbf{16} (1977), 423-433

\bibitem {Brieskorn} E. Brieskorn, \textit{Rationale Singularit\"aten komplexer Fl\"achen},
 Invent. Math. \textbf{4} (1968), 336-358


\bibitem {KM} S. Keel and J. McKernan, \textit{Rational curves on quasi-projective surfaces},
Mem. Amer. Math. Soc. \textbf{140} (1999), no. 669

\bibitem {Keum} J. Keum, \textit{A rationality criterion for projective surfaces - partial solution to Kollar's conjecture},
Algebraic geometry, 75-87, Contemp. Math. \textbf{422}, Amer.
Math. Soc., Providence, RI, 2007

\bibitem {KZ} J. Keum, and D.-Q. Zhang, \textit{Fundamental groups of open K3 surfaces, Enriques surfaces and Fano 3-folds},
Jour. Pure Applied Algebra \textbf{170} (2002), 67-91

\bibitem {KNS} R. Kobayashi, S. Nakamura, and F. Sakai \textit{A
numerical characterization of ball quotients for normal surfaces
with branch loci}, Proc. Japan Acad. Ser. A,  Math. Sci.
\textbf{65} (1989), no. 7, 238-241

\bibitem {Kollar05} J. Koll\'ar, \textit{Einstein metrics on $5$--dimensional Seifert bundles}, Jour. Geom. Anal.
\textbf{15} (2005), no. 3, 463-495

\bibitem {Kollar06} J. Koll\'ar, \textit{Is there a topological Bogomolv-Miyaoka-Yau inequality?} to
appear, math.AG/0602562

\bibitem {KSB} J. Koll\'ar, and N. I. Shepherd-Barron, \textit{Threefolds and deformations of
surface singularities} Invent. Math. \textbf{91} (1988), no. 2,
299--338

\bibitem {Kon} S. Kond\=o, \textit{Enriques surfaces with finite automorphism groups},
Japanese J. Math. \textbf{12} (1986), 191-282

\bibitem {LW} E. Looijenga, and J. Wahl, \textit{Quadratic functions and smoothing surface singularities},
Topology \textbf{25} (1986), no.3, 261-291

\bibitem {Manetti} M. Manetti, \textit{Normal degenerations of the complex projective plane},
J. Reine Angew. Math. \textbf{419} (1991), 89-118

\bibitem {Matsuki} K. Matsuki, \textit{Introduction to the Mori program},
Universitext. Springer-Verlag, New York, 2002

\bibitem {Megyesi} G. Megyesi, \textit{Generalisation of the Bogomolov-Miyaoka-Yau inequality to singular surfaces},
Proc. London Math. Soc. (2) \textbf{78} (1999), 241-282


\bibitem {Miyaoka} Y. Miyaoka, \textit{The maximal number of quotient singularities on surfaces with given numerical invariants},
Math. Ann. \textbf{268} (1984), 159-171



\bibitem {MY} D. Montgomery and C. T. Yang, \textit{Differentiable pseudo-free circle
actions on homotopy seven spheres}, Proc. of the Second Conference
on Compact Transformation Groups (Univ. Massachusetts, Amherst,
Mass., 1971), Part I, Springer, 1972, pp. 41-101. Lecture Notes in
Math., Vol. 298

\bibitem {N} V. V. Nikulin, \textit{Integral symmetric bilinear forms and its applications}, Izv. Akad. Nauk SSSR Ser. Mat. \textbf{43} (1979), no. 1,
111-177; English translation: Math. USSR Izv. \textbf{14} (1979),
no. 1, 103-167 (1980)

\bibitem {Rie77} O. Riemenschneider,  \textit{Die Invarianten der endlichen Untergruppen von GL(2, $\mathbb{C}$)},
Math. Zeit. \textbf{153} (1977), 37-50

\bibitem {Sakai} F. Sakai, \textit{Semistable curves on algebraic surfaces and logarithmic pluricanonical
maps}, Math. Ann. \textbf{254} (1980), no. 2, 89-120
\bibitem {Serre} J. P. Serre, \textit{A Course in Arithmetic},
 Springer-Verlag, New York, 1973

\end{thebibliography}

\end{document}